\documentclass[12pt]{article}
\usepackage{amsfonts}
\usepackage{amsthm}
\usepackage{amsmath}
\usepackage{mathrsfs}
\usepackage{amssymb}
\usepackage{graphicx}
\usepackage{caption}
\RequirePackage[colorlinks,linkcolor=red,citecolor=red,urlcolor=red]{hyperref}
\usepackage{latexsym,epsf,dsfont,color,eurosym,mathrsfs,url}

\newlength{\fixboxwidth}
\newtheorem{lemma}{Lemma}
\newtheorem{theorem}{Theorem}

\newcommand{\RR}{{\mathbb R}}
\newcommand{\EE}{{\mathbb E}}
\newcommand{\II}{{\mathbb I}}

\newcommand{\PP}{\mathbb P}
\newcommand{\bfP}{\mathbf P}

\newcommand{\argmin}[1]{\underset{#1}{\mathrm{argmin}}}

\def\Acal{\mathcal A}
\def\Bcal{\mathcal B}

\def\Ecal{\mathcal E}
\def\Gcal{\mathcal G}
\def\Hcal{\mathcal H}
\def\Kcal{\mathcal K}
\def\Ncal{\mathcal N}
\def\Rcal{\mathcal R}

\def\Wcal{\mathcal W}
\def\Xcal{\mathcal X}

\def\hat{\widehat}
\def\epsilon{\varepsilon}
\def\argmin{\mathop{\rm argmin}}

\begin{document}

\title{Minimax Optimal Rates of Estimation in High Dimensional Additive Models: Universal Phase Transition}

\date{(\today)}

\author{\\
Ming Yuan$^\ast$ and Ding-Xuan Zhou$^\dag$\\
\\
University of Wisconsin-Madison and City University of Hong Kong}

\footnotetext[1]{
Department of Statistics, University of Wisconsin-Madison, 1300 University Avenue, Madison, WI 53706. The research of Ming Yuan was supported in part by NSF Career Award DMS-1321692, FRG Grant DMS-1265202 and NIH Grant 1-U54AI117924-01.}
\footnotetext[2]{
Department of Mathematics, City University of Hong Kong, Tat Chee Avenue, Kowloon, Hong Kong. The research of Ding-Xuan Zhou was supported in part by the Research Grants Council of Hong Kong [Project No. CityU 104012] and by National Natural Science Foundation of China [Project No. 11461161006].}
\maketitle

\newpage

\begin{abstract}
We establish minimax optimal rates of convergence for estimation in a high dimensional additive model assuming that it is approximately sparse. Our results reveal an interesting phase transition behavior universal to this class of high dimensional problems. In the {\it sparse regime} when the components are sufficiently smooth or the dimensionality is sufficiently large, the optimal rates are identical to those for high dimensional linear regression, and therefore there is no additional cost to entertain a nonparametric model. Otherwise, in the so-called {\it smooth regime}, the rates coincide with the optimal rates for estimating a univariate function, and therefore they are immune to the ``curse of dimensionality''.
\end{abstract}


\newpage

\section{Introduction}
With the recent advances in science and technology, high dimensional regression problems have become ubiquitous in a multitude of areas -- genomics, medical imaging, and finance are a few well known examples. Considerable amount of research effort has been devoted to the understanding of challenges brought about by the high dimensionality, and development of statistical methodology to counter them. Most of the existing work focuses on high dimensional linear regression where a number of approaches such as the nonnegative garrote (Breiman, 1995), the Lasso (Tibshirani, 1996), the SCAD (Fan and Li, 2001), and the Dantzig selector (Cand\`es and Tao, 2007), have been developed to exploit sparsity, or perform variable selection; and much progress has also been made to understand to what extent a high dimensional regression coefficient vector can be reliably estimated; see, e.g., Koltchinskii (2011), B\"uhlmann and van de Geer (2013) and references therein.

Linear models, however, could be too restrictive in many applications. As a more flexible alternative, high dimensional additive models have attracted much attention in the past several years. See, e.g., Lin and Zhang (2006), Yuan (2007), Koltchinskii and Yuan (2008), Ravikumar et al. (2008), Meier, van de Geer and B\"uhlmann (2009), Koltchinskii and Yuan (2010) and Raskutti, Wainwright and Yu (2012) among others. Let $\{(X_i,Y_i): i=1,\dots ,n\}$ be independent copies of a random couple $(X,Y)$ following a regression model:
\begin{equation}
\label{eq:regress}
Y=f(X)+\epsilon,
\end{equation}
where the error $\epsilon$ follows $\Ncal(0,\sigma^2)$ distribution. The additive model amounts to the assumption that
\begin{equation}
\label{eq:add}
f(x_1,\ldots,x_d)=f_1(x_1)+\cdots+f_d(x_d),
\end{equation}
where the component functions $f_j$s are modeled non-parametrically; see, e.g., Stone (1985) or Hastie and Tibshirani (1990). Here we assume that they reside in certain reproducing kernel Hilbert spaces (RKHS); see, e.g., Aronszajn (1950) and Wahba (1990).

To fix ideas, assume that $X$ follows a distribution $\Pi$ supported on a product space $\Xcal^d$ for some compact subset $\Xcal$ of $\RR$; and that all component functions come from a common RKHS of functions on $\Xcal$, denoted by $(\Hcal_1, \|\cdot\|_{\Hcal_1})$. It is clear that the additive model (\ref{eq:add}) can be identified with space
\begin{eqnarray*}
\Hcal_d := \Hcal_1 \oplus \ldots \oplus \Hcal_1 =\biggl\{g:\Xcal^d\to \RR| g(x_1,\ldots,x_d)=g_1 (x_1) + \ldots + g_d (x_d), \\
{\rm\ and\ } g_1, \ldots, g_d \in \Hcal_1\biggr\}.
\end{eqnarray*}
Obviously linear models can be viewed as a trivial special case of (\ref{eq:add}) by taking $\Hcal_1$ to be the collection of all univariate linear functions defined over $\Xcal$. Another canonical example of $\Hcal_1$ is the $\alpha$th ($\alpha>1/2$) order Sobolev space $\Wcal_2^\alpha([0,1])$ defined on a unit interval ($\Xcal=[0,1]$). See, e.g., Wahba (1990) for further examples.

We note that for a general $g\in \Hcal_d$, the additive representation given by (\ref{eq:add}) may not be unique. Define the (quasi-)norm $\|f\|_{\ell_q (\Hcal_d)}$ ($q> 0$) by
\begin{eqnarray*}
\|g\|_{\ell_q (\Hcal_d)} = \inf \biggl\{\left\|(\|g_1\|_{\Hcal_1},\ldots,\|g_d\|_{\Hcal_1})^\top\right\|_{\ell_q}: g_1 (x_1) + \ldots + g_d (x_d) = g(x_1, \ldots, x_d)\\
{\rm\ and\ } g_1, \ldots, g_d \in \Hcal_1\biggr\}.
\end{eqnarray*}
In other words, $\|f\|_{\ell_q (\Hcal_d)}$ is the $\ell_q$ norm of the vector of RKHS norms of its component functions minimized over all of its additive representations. In particular, when $q=2$, $\|\cdot\|_{\ell_2 (\Hcal_d)}$ can be viewed as a RKHS norm. More specifically, let $K: \Xcal \times \Xcal \to \RR$ be a Mercer kernel generating the RKHS $(\Hcal_1, \|\cdot\|_{\Hcal_1})$ and write
$$
K_d((x_1, \ldots, x_d)^\top,(x_1', \ldots, x_d')^\top)=K(x_1,x_1')+\cdots+K(x_d,x_d').
$$
It is not hard to see that $K_d$ is the generating kernel of RKHS $(\Hcal_d,\|\cdot\|_{\ell_2 (\Hcal_d)})$. Another special case of the $\ell_q(\Hcal_d)$ norm defined above is the case when $q\downarrow 0$. $\|\cdot\|_{\ell_0(\Hcal_d)}$ can be interpreted as the smallest number of additive components needed to express a function from $\Hcal_d$.

\begin{figure}[htbp]
\begin{center}
\includegraphics[width=0.8\textwidth]{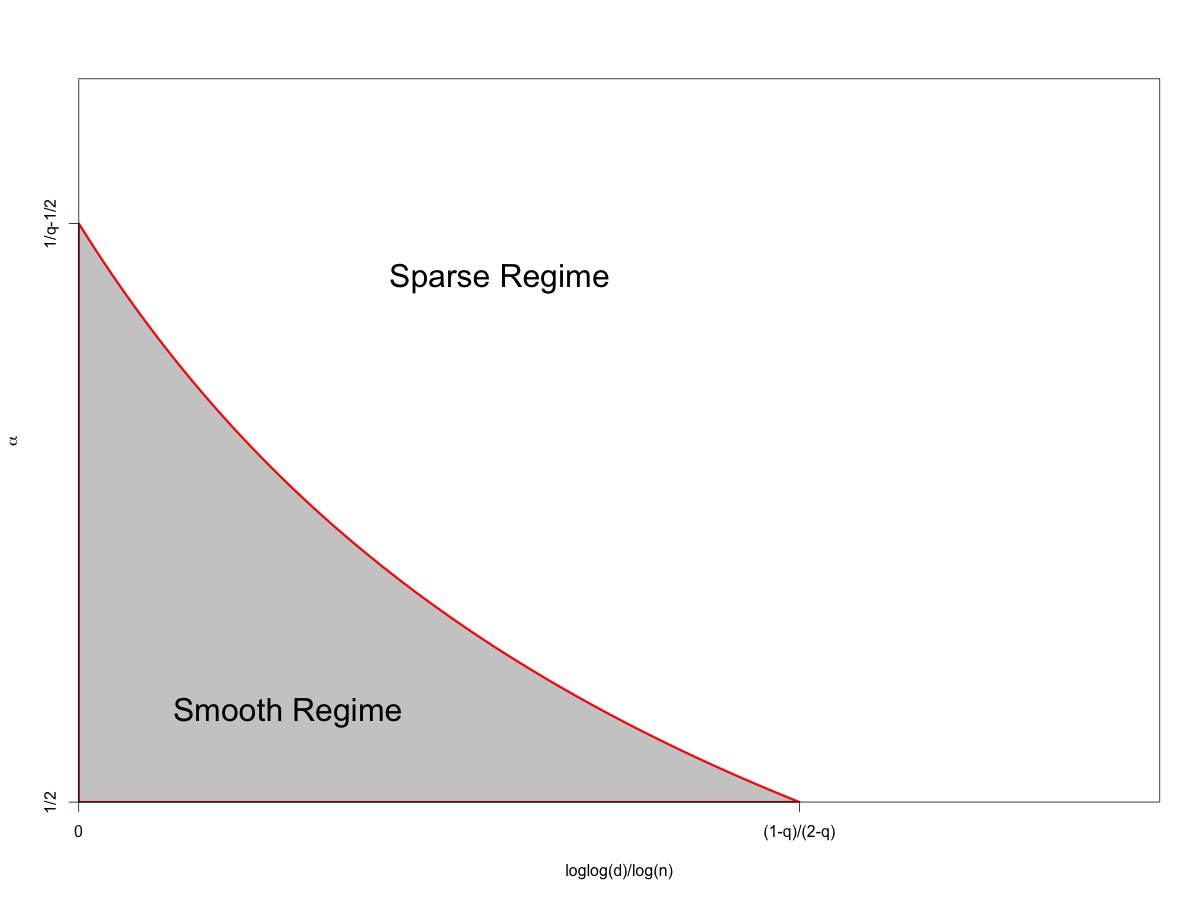}
\caption{Phase transition in optimal rates of convergence: When the smoothness index $\alpha$ and dimensionality measured by $\log \log d/\log n$ falls in the {\it smooth region} in the figure above, the optimal rate is given by $n^{-{2\alpha\over 2\alpha+1}}$ which is determined solely by the smoothness index. On the other hand, if they fall into the {\it sparse regime}, then the optimal rate is given by $(n^{-1}\log d)^{1-q/2}$ which is determined entirely by the dimensionality.}
\label{fig:rate}
\end{center}
\end{figure}

When the dimension $d$ is large, it is of particular interests to consider the case when $f$ resides in an $\ell_q(\Hcal_d)$ ball for $0< q< 1$:
$$
\Bcal_R\left(\ell_q (\Hcal_d)\right)=\left\{g\in \Hcal_d:\, \|g\|_{\ell_q(\Hcal_d)}^q\le R\right\}.
$$
Write
$$
\|g\|_{L_2(\Pi)}=\left(\int_{\Xcal^d} g^2(x)d\Pi(x)\right)^{1/2}.
$$
We are interested in the minimax optimal rate of convergence for estimating $f$ in terms of the squared $\|\cdot\|_{L_2(\Pi)}$ norm. In particular, when $\Hcal_1$ is taken to be the $\alpha$th order Sobolev space $\Wcal_2^\alpha$ defined on the unit interval, our results imply that the minimax optimal rate for estimating $f\in \Bcal_R(\ell_q(\Hcal_d))$ is given by
\begin{equation}
\label{eq:optrate}
\Rcal(n,d)=\left({\log d\over n}\right)^{1-{q\over 2}}+n^{-{2\alpha\over 2\alpha+1}},
\end{equation}
up to a constant scaling factor. The optimal rate of convergence given by (\ref{eq:optrate}) exhibits an interesting phase transition phenomenon as illustrated in Figure \ref{fig:rate}. 

More specifically, when the component functions are not sufficiently smooth in that
$$
\alpha<{1\over q}-{1\over 2},
$$
the second term on the right hand side of (\ref{eq:optrate}) is dominated by the first one if $d$ is ultra-large:
$$
d> \exp\left[n^{{2\over 2-q}\left({1\over 2\alpha+1}-{q\over 2}\right)}\right],
$$
and hence the minimax optimal rate becomes
\begin{equation}
\label{eq:optratepar}
\Rcal(n,d)\asymp\left({\log d\over n}\right)^{1-{q\over 2}},
\end{equation}
where we write for two positive sequences $a_{n,d}$ and $b_{n,d}$, $a_{n,d}\asymp b_{n,d}$ if $a_{n,d}/b_{n,d}$ is bounded away from both zero and infinity. The rate given by (\ref{eq:optratepar}) happens to be the minimax optimal optimal rate for estimating a $d$ dimensional linear regression when assuming the vector of regression coefficient comes from a $\ell_q$ ball in $\RR^d$; see, e.g., Ye and Zhang (2010) or Raskutti, Wainwright and Yu (2011). On the other hand, when
$$
d\le \exp\left[n^{{2\over 2-q}\left({1\over 2\alpha+1}-{q\over 2}\right)}\right],
$$
the optimal rate is given by
$$
\Rcal(n,d)\asymp n^{-{2\alpha\over 2\alpha+1}}.
$$
This rate coincides with the optimal rate for estimating $f$ if we know in advance that it actually comes from a single component space $\Hcal_1$, e.g., $f_2=\cdots=f_d=0$, rather than the $d$-variate function space $\Hcal_d$; see, e.g., Stone (1980, 1982) and Tsybakov (2009). Similar phase transition depending on the dimensionality $d$ has also been observed earlier for high dimensional additive models under exact sparsity ($q=0$). See, e.g., Koltchinskii and Yuan (2010), Raskutti, Wainwright and Yu (2012) and Suzuki and Sugiyama (2013). Our results suggest that such phenomenon is more universal and applies in general to the approximate sparse case.

It is also worth pointing out that such a phase transition in $d$ vanishes when the component functions are sufficiently smooth in that
$$
\alpha\ge {1\over q}-{1\over 2},
$$
a phenomenon absent in the case of exact sparsity ($q=0$). In this situation, the second term on the right hand side of (\ref{eq:optrate}) is always dominated by the first one and therefore the optimal rate is always
$$
\Rcal(n,d)\asymp\left({\log d\over n}\right)^{1-{q\over 2}}.
$$
In other words, we pay no extra price, in terms of rates of convergence, for entertaining a generally nonparametric additive model (\ref{eq:add}) when compared with the much more restrictive linear models, regardless of the value of $d$.

Although we focus on additive models, our general framework is also closely related to multiple kernel learning or ``aggregation" of kernel machines, a popular technique in machine learning to combine multiple kernels instead of using a single one in order to achieve improved prediction performance. These type of problems have been studied recently by Bousquet et al. (2003), Cramer et al. (2003), Lanckriet et al. (2004), Micchelli and Pontil (2005), Srebro and Ben-David (2006), Bach (2008), and Suzuki and Sugiyama (2013) among others. It is expected that our results here could lead to further understanding of these problems as well.

The rest of the paper is organized as follows. We first review some basic concepts and properties of reproducing kernel Hilbert spaces in Section \ref{sec:prep}. Section \ref{sec:main} presents the main results. All proofs are relegated to Section \ref{sec:proof}.

\section{Reproducing Kernel Hilbert Spaces}
\label{sec:prep}

We begin with a brief review of some of the basic facts about RKHS which we shall make repeated use later on. Interested readers are referred to Aronszajn (1950) and Wahba (1990) for further details. In particular, we shall focus on the $j$th component space, e.g., the RKHS defined on the $j$th coordinate of $X\in \Xcal^d$.

\subsection{Kernel and RKHS}
Recall that $K$ is a symmetric positive semi-definite, square integrable function on $\Xcal\times \Xcal$. It can be uniquely identified with the Hilbert space $\Hcal_1$ that is the completion of 
$$
\{K(x,\cdot): x\in \Xcal\}
$$
under the inner product
$$
\left\langle \sum_i c_i K(x_i,\cdot),\sum_j c_j' K(x_{j}',\cdot)\right\rangle_K =\sum_{i,j}c_ic_j'K(x_i,x_j').
$$
In the rest of the section, we shall write $\Hcal_1$ and $\Hcal(K)$ interchangeably with the latter notion emphasizing the one-to-one correspondence between a kernel and a RKHS. Most, if not all, the commonly used kernels are bounded, which we shall assume in what follows. In fact, without loss of generality, we shall assume in the rest of the paper that $\sup_x K(x,x)=1$. Note that, for any $h\in \Hcal(K)$,
\begin{equation}
\label{eq:definf}
\|h\|_\infty:=\sup_{x\in \Xcal}|h(x)|=\sup_{x\in \Xcal}\left|\langle h, K(x,\cdot)\rangle_K\right|\le \sup_x\|K(x,\cdot)\|_K\|h\|_K,
\end{equation}
by Cauchy-Schwartz inequality. Recall that
$$
\|K(x,\cdot)\|_K^2=\langle K(x,\cdot),K(x,\cdot)\rangle_K=K(x,x)\le 1.
$$
Thus,
$$
\|h\|_\infty\le\|h\|_K,
$$
a convenient fact that we shall used repeatedly in the later analysis.

By spectral theorems, $K$ admits the following eigenvalue decomposition:
$$
K(x,x')=\sum_{k\ge 1}\lambda_{jk} \varphi_{jk}(x)\varphi_{jk}(x')
$$
where $\lambda_{j1}\ge\lambda_{j2}\ge\cdots\ge 0$ are its eigenvalues and $\{\varphi_{jk}: k\ge 1\}$ are the corresponding eigenfunctions such that
$$
\langle \varphi_{jk},\varphi_{jk'}\rangle_{L_2(\Pi_j)}=\delta_{kk'}.
$$
Here $\Pi_j$ is the $j$th marginal distribution of $\Pi$, and $\delta_{kk'}$ is the Kronecker's delta. It is well known that the RKHS-norm of any $h\in \Hcal(K)$ can be written as 
$$
\|h\|_{K}^2 = \sum_{k\geq 1} {1\over {\lambda_{jk}}}\langle h, \varphi_{jk}\rangle_{L_2(\Pi_j)}^2, 
$$
which means that the ``smoothness'' of functions in $\Hcal(K)$ are determined by the rate of decay of eigenvalues $\lambda_{jk}$, and the unit balls in the RKHS $\Hcal(K)$ are ellipsoids in the space $L_2(\Pi_j)$ with ``axes'' $\sqrt{\lambda_{jk}}$. For example, it is well known that if $\Pi_j$ is the Lebesgue measure on $[0,1]$, then $\lambda_{jk}\asymp k^{-2\alpha}$ for $\Wcal_2^\alpha$.

\subsection{Complexity of RKHS}
How well we can recover a function from a particular RKHS is fundamentally related to the capacity of the unit ball in $\Hcal(K)$:
$$
\Bcal_1(\Hcal(K)):=\left\{h\in \Hcal(K): \|h\|_K\le 1\right\}.
$$
See, e.g., Yang and Barron (1999). In particular, the capacity of $\Bcal_1(\Hcal(K))$ can be measured by its covering number $\Ncal(\Bcal_1(\Hcal(K)), \delta, \|\cdot\|_{\infty})$ where $\|\cdot\|_\infty$ is defined in (\ref{eq:definf}). Recall that for $\delta>0$ and a set ${\mathcal F}$ of continuous functions on a metric space $\Xcal$, the covering number ${\mathcal N} ({\mathcal F}, \delta, \|\cdot\|_{\infty})$ with respect to the $\|\cdot\|_\infty$ metric is defined as the smallest integer $m$ such that
$${\mathcal F} = \bigcup_{i=1}^m \{f\in {\mathcal F}: \|f- f^{(i)}\|_{\infty} \le \delta\}$$
for some $\{f^{(i)}\}_{i=1}^m \subset {\mathcal F}$. In particular, if $\lambda_{jk}=O(k^{-2\alpha})$ and $\sup_{k\ge 1}\|\varphi_{jk}\|_{\infty}<\infty$, then
\begin{equation}
\label{eq:capassume}
\log \Ncal(\Bcal_1(\Hcal(K)), \delta, \|\cdot\|_{\infty}) \leq c \delta^{-{1\over \alpha}}, \qquad \forall \delta >0,
\end{equation}
for some constant $c>0$. This holds, for example, for Sobolev spaces of order $\alpha$.

For our purposes, we are also interested in certain data-dependent estimates of the complexity of a function class, namely, Rademacher and Gaussian complexities. See, e.g., Bartlett and Mendelson (2002). Write
\begin{equation}
\label{eq:radcomp}
R_{jn}(u):=\sup_{h\in \Bcal_1(\Hcal(K)): \|h\|_{L_2(\Pi_{j})}\le u}\left|{1\over n}\sum_{i=1}^n \sigma_i h(x_{ij})\right|
\end{equation}
where $\sigma_i$s are iid Rademacher variables, that is $\PP(\sigma_i=1)=\PP(\sigma_i=-1)=1/2$. The following bound of $R_{jn}$ will become useful for our later analysis.

\begin{lemma}
\label{le:emp1}
Assume that $\lambda_{jk}\le c_1k^{-2\alpha}$ and $\sup_{k\ge 1}\|\varphi_{jk}\|_{L_\infty}<c_2$ for some constants $c_1,c_2>0$. Then there exists a constant $c>0$ depending on $\alpha$, $c_1$ and $c_2$ only such that for any $\beta>0$, with probability at least $1-d^{-\beta}$,
$$
R_{jn}(u)\le cn^{-1/2}\left(u^{1-{1\over 2\alpha}}+u\sqrt{\beta\log d}+{\beta\log d\over \sqrt{n}}+e^{-d}\right)
$$
uniformly for all $u\in [0,1]$.
\end{lemma}

Another quantity of interests to us is the ``empirical'' Gaussian complexity of the unit ball in $\Hcal(K)$:
\begin{equation}
\label{eq:gausscomp}
\hat{Z}_{jn}(u):=\sup_{h\in \Bcal_1(\Hcal(K)): \|h\|_{L_2(\Pi_{jn})}\le u}\left|{1\over n}\sum_{i=1}^n \epsilon_i h(x_{ij})\right|
\end{equation}
where $\Pi_{jn}$ is the $j$th marginal of the empirical distribution $\Pi_n$. Similar to Lemma \ref{le:emp1}, we have the following bound for $\hat{Z}_{jn}$.

\begin{lemma}
\label{le:emp}
Assume that $\lambda_{jk}\le c_1k^{-2\alpha}$ and $\sup_{k\ge 1}\|\varphi_{jk}\|_{L_\infty}<c_2$ for some constants $c_1,c_2>0$. Then there exists a constant $c>0$ depending on $\alpha$, $c_1$ and $c_2$ only such that for any $\beta>0$, with probability at least $1-d^{-\beta}$,
$$
\hat{Z}_{jn}(u)\le cn^{-1/2}\left(u^{1-{1\over 2\alpha}}+u\sqrt{\beta\log d}+e^{-d}\right)
$$
uniformly for all $u\in [0,1]$.
\end{lemma}

Both Lemmas \ref{le:emp1} and \ref{le:emp} follow from a standard peeling argument (see, e.g., van de Geer, 2000). We present their proofs in the Appendix for completeness.

\section{Main Results}
\label{sec:main}
In what follows, we shall assume that there exists a constant $\eta_q>1$ such that
\begin{equation}
\label{eq:reassumption}
\eta_q^{-1}\|g\|_{L_2(\Pi)}^2\le \sum_{j=1}^d\|g_j\|_{L_2(\Pi_j)}^2\le \eta_q\|g\|_{L_2(\Pi)}^2
\end{equation}
for any $g\in \Bcal_R(\ell_q(\Hcal_d))$, where
$$
g(x_1,\ldots,x_d)=g_1(x_1)+\cdots+g_d(x_d)
$$
and
$$
\|g\|_{\ell_q(\Hcal_d)}^q=\sum_{j=1}^d \|g_j\|_{\Hcal_1}^q.
$$
Condition (\ref{eq:reassumption}) is a nonparametric version of the restricted eigenvalue conditions commonly used in analyzing sparse estimation in high dimensional linear regression; see, e.g., Bickel, Ritov and Tsybakov (2009). It is worth noting that different from the usual restricted eigenvalue conditions in linear regression, Condition (\ref{eq:reassumption}) is on the {\it distribution} of $X$ rather than the design matrix, or observations $X_1,\ldots,X_n$. The condition is satisfied in particular when $\Pi$ is a product measure.

To fix ideas, in the rest of the paper, we shall also assume that there exist a constant $c_\lambda>1$ and a non-increasing sequence of nonnegative numbers $\lambda_1\ge\lambda_2\ge \cdots$ such that
\begin{equation}
\label{eq:eigassumption}
c_\lambda^{-1}\lambda_k\le \lambda_{jk}\le c_\lambda\lambda_k,
\end{equation}
for all $j=1,2,\ldots, d$ and $k\ge 1$. In addition, similar to the treatment of high dimensional linear models (see, e.g., Raskutti, Wainwright and Yu, 2011), we shall assume in the rest of the paper that $c_0n^{q/2}\le d\le e^n$ for some universal constant $c_0>0$ to ensure nontrivial probabilistic bounds. This, in particular, is true in high dimensional settings where $n<d<e^n$.

We are now in position to present the main results. We first state a minimax lower bound. 
\begin{theorem}
\label{th:lower}
Assume that $\lambda_k=k^{-2\alpha}$ for some $\alpha>1/2$. Under the regression model (\ref{eq:regress}) where $f\in \Bcal_R(\ell_q(\Hcal_d))$ and the covariate $X$ follows a distribution $\Pi$ such that (\ref{eq:reassumption}) and (\ref{eq:eigassumption}) hold,  and the eigenfunctions $\{\varphi_{jk}: j=1,\ldots, d, k\ge 1\}$ are uniformly bounded, there exists a constant $c>0$ depending on $\sigma^2$, $\alpha$, $R$, $c_\lambda$ and $\eta_q$ only such that
$$
\lim_{n\to\infty}\inf_{\tilde{f}}\sup_{f\in \Bcal_R(\ell_q(\Hcal_d))}\PP\left\{\|\tilde{f}-f\|_{L_2(\Pi)}^2\ge c\left[\left(\log d\over n\right)^{1-q/2}+n^{-{2\alpha\over 2\alpha+1}}\right]\right\}>0.
$$
\end{theorem}

The lower bound is established via Fano's Lemma. See, e.g., Cover and Thomas (1991). We relegate its proof to Section \ref{sec:proof}. Next, we show that the rates given in the lower bound in the previous theorem is attainable. In particular, we consider the least squares estimator:
\begin{equation}
\label{eq:lse}
\hat{f}=\argmin_{g\in \Bcal_R(\ell_q(\Hcal_d))} \left\{{1\over n}\sum_{i=1}^n \left[Y_i-g(X_i)\right]^2\right\}.
\end{equation}
The next result shows that $\hat{f}$ is indeed minimax rate optimal.

\begin{theorem}
\label{th:upper}
Assume that $\lambda_k=k^{-2\alpha}$ for some $\alpha>1/2$. Under the regression model (\ref{eq:regress}) where $f\in \Bcal_R(\ell_q(\Hcal_d))$ and the covariate $X$ follows a distribution $\Pi$ such that (\ref{eq:reassumption}) and (\ref{eq:eigassumption}) hold,  and the eigenfunctions $\{\varphi_{jk}: j=1,\ldots, d, k\ge 1\}$ are uniformly bounded, there exists a constant $c>0$ depending on $\sigma^2$, $\alpha$, $R$, $c_\lambda$ and $\eta_q$ only such that for any $\beta>0$ with probability at least $1-d^{-\beta}$,
\begin{equation}
\label{eq:bdlse}
\|\hat{f}-f\|_{L_2(\Pi)}^2\le c(\beta+1)\left[\left(\log d\over n\right)^{1-q/2}+n^{-{2\alpha\over 2\alpha+1}}\right],
\end{equation}
and
\begin{equation}
\label{eq:bdlse2}
\|\hat{f}-f\|_{L_2(\Pi_n)}^2\le c(\beta+1)\left[\left(\log d\over n\right)^{1-q/2}+n^{-{2\alpha\over 2\alpha+1}}\right],
\end{equation}
where $\hat{f}$ is the least squares estimator defined by (\ref{eq:lse}).
\end{theorem}

The proof of Theorem \ref{th:upper} is also presented in Section \ref{sec:proof}. It relies on several basic facts of the empirical processes theory such as symmetrization inequalities and contraction inequalities for Rademacher processes that can be found in the books of Ledoux and Talagrand (1991) and van der Vaart and Wellner (1996). We also use TalagrandÕs concentration inequality for empirical processes; see, e.g., Talagrand (1996) and Bousquet (2002).

Theorems \ref{th:lower} and \ref{th:upper} together immediate imply that the minimax optimal rate for estimating $f\in \Bcal_R(\ell_q(\Hcal_d))$ is
$$
\|\hat{f}-f\|_{L_2(\Pi)}^2\asymp \left(\log d\over n\right)^{1-q/2}+n^{-{2\alpha\over 2\alpha+1}}.
$$
This result connects with two strands of literature -- estimating high dimensional linear regression assuming that the coefficient vector belongs to an $\ell_q$ ball, and estimating a high dimensional additive model assuming that the underlying function comes from a $\ell_0(\Hcal_d)$ ball. In the case of linear regression, it is known that $\ell_1$ penalty or the Lasso (Tibshirani, 1996) leads to rate optimal estimators under suitable regularity conditions. See, e.g., Ye and Zhang (2010). Similar phenomenon has also been observed for the high dimensional additive models where it is shown that a mixed $\ell_1$ norm penalty of the form
\begin{equation}
\label{eq:doublepen}
a_n^2\sum_{j=1}^d \|g_j\|_{\Hcal_1}+a_n\sum_{j=1}^d \|g_j\|_{L_2(\Pi_{jn})}
\end{equation}
can lead to rate optimal estimators with appropriate choices of the tuning parameter $a_n>0$. See, e.g., Koltchinskii and Yuan (2010) and Raskutti, Wainwright and Yu (2012). The use of a mixed $\ell_1$ penalty of the form (\ref{eq:doublepen}) highlights the difference between linear models and additive models. When dealing with nonparametric component functions, we need to penalize both the RKHS norm and $L_2$ norm, the former ensures smoothness of the estimate whereas the latter is needed for thresholding redundant components and hence sparsity.

A natural question is whether or not a similar strategy will lead to minimax rate optimal estimators under $\ell_q(\Hcal_d)$ ball for general $0<q<1$. Somewhat surprisingly, the answer appears to be negative in general. And we give here a heuristic argument why. The challenge occurs in the smooth regime where
$$
\alpha<{1\over q}-{1\over 2},\qquad {\rm and}\qquad d\le \exp\left[n^{{2\over 2-q}\left({1\over 2\alpha+1}-{q\over 2}\right)}\right].
$$
Recall that the corresponding minimax optimal rate of convergence in the smooth regime is given by
$$
n^{-{2\alpha\over 2\alpha+1}}.
$$
As pointed out before, this is the best possible rate of convergence even if there is only one nonzero component. And to achieve this rate, we need to choose 
\begin{equation}
\label{eq:tune}
a_n\gtrsim n^{-{\alpha\over 2\alpha+1}},
\end{equation}
because, if $a_n$ is smaller, then in the particular case of one nonzero component, the minimax optimal rate cannot be attained. See, e.g., Tsybakov (2009) or Koltchinskii and Yuan (2010).  Now for a general $f$ from the unit $\ell_q(\Hcal_d)$ ball, we will need a diverging number of nonzero components to approximate it. More precisely, as we shall show in the proofs, we may need estimate up to
$$
\left\lceil\left(n\over \log d\right)^{q/2}\right\rceil
$$
nonzero components to balance the approximation error and sparsity. If we choose $a_n$ to be of the order given by (\ref{eq:tune}), then each component can only be estimated with squared $L_2$ error of the order of
$$a_n^2\gtrsim n^{-{2\alpha\over 2\alpha+1}},$$
leading to an overall rate of convergence no smaller than, up to a multiplicative constant,
$$
\left(n\over \log d\right)^{q/2}n^{-{2\alpha\over 2\alpha+1}},
$$
at least under the assumption that $\Pi$ is a product measure. This rate is obviously suboptimal. As a result, in the smooth regime, no matter what value $a_n$ is, we cannot attain the minimax optimal rate of convergence through a mixed $\ell_1$ penalty of the form (\ref{eq:doublepen}).

\section{Proofs.}
\label{sec:proof}
We now prove the main results Theorems \ref{th:lower} and \ref{th:upper}. For brevity, we shall also assume that $\sigma^2=1$ and $R=1$ in the proofs. The more general case follows an identical arguments with different constants.

\subsection{Lower bounds.}
We establish the lower bound via Fano's Lemma.  To this end, we need to construct a set of functions
$$\Gcal:=\{g^1,\ldots, g^M\}\subset \Bcal_1(\ell_q(\Hcal_d))$$
that are sufficiently apart from each other. Let $N$ be a natural number whose value will be specified later. For a matrix $A\in\{-1,0,1\}^{d\times N}$, denote by $s_A$ the number of its nonzero rows, that is
$$
s_A={\rm card}\left\{i: A_{i\cdot}\neq {\bf 0}\right\},
$$
where $A_{i\cdot}$ is the $i$th row vector of $A$. Write
$$
g_A(x_1,\ldots,x_d)=N^{-1/2}s_A^{-1/q}\sum_{j=1}^d\sum_{k=1}^N a_{jk}\lambda^{1/2}_{j,N+k}\varphi_{j,N+k}(x_j).
$$
It is clear that
\begin{eqnarray*}
\|g_A\|_{\ell_q(\Hcal_d)}^q&\le& N^{-q/2}s_A^{-1}\sum_{j=1}^d\left\|\sum_{k=1}^N a_{jk}\lambda^{1/2}_{j,N+k}\varphi_{j,N+k}(x_j)\right\|_{\Hcal_1}^q\\
&=&s_A^{-1}\sum_{j=1}^d \left(N^{-1}\sum_{k=1}^N a_{jk}^2\right)^{q/2}.
\end{eqnarray*}
Because $a_{jk}^2\in \{0,1\}$, this can be further bounded by
$$
\|g_A\|_{\ell_q(\Hcal_d)}^q\le s_A^{-1}\sum_{j=1}^d \II(A_{i\cdot}\neq 0)=1,
$$
which implies that $g_{A}\in \Bcal_1(\ell_q(\Hcal_d))$.

We now describe how to generate the set $\Gcal$. In particular, we consider functions of the form $g_A$ with $A\in \{\pm 1, 0\}^{d\times N}$ as described before. We first choose $s$ rows of $A$ to be nonzero, and set the rest of the rows of $A$ to be zero. The value of $s$ will become clear later. To this end, we appeal to Vershamov-Gilbert Lemma which states that we can find a set $\{\theta_1,\ldots,\theta_{M_1}\}\subset\{0,1\}^d$ such that
\begin{enumerate}
\item[(a)] $\|\theta_k\|_{\ell_1}=s$ for $1\le k\le M_1$;
\item[(b)] for any $k\neq k'$, $\|\theta_k-\theta_{k'}\|_{\ell_1}\ge s/2$;
\item[(c)] $\log M_1\ge {1\over 4}s\log (d/s)$.
\end{enumerate}
See, e.g., Massart (2007). For a given $\theta$, we set zero the rows of $A$ if the corresponding coordinate of $\theta$ is zero. In the next step, we fill in the remaining rows of $A$ with $\pm 1$. Again, by Vershamov-Gilbert Lemma, there exists a set $\{\Gamma_1,\ldots,\Gamma_{M_2}\}\in \{\pm 1\}^{s\times N}$ such that
\begin{enumerate}
\item[(a')] for any $k\neq k'$, $\|\Gamma_k-\Gamma_{k'}\|_{\rm F}^2\ge Ns/2$;
\item[(b')] $\log M_2\ge Ns/8$.
\end{enumerate}
For a given $\Gamma$, we shall fill in the nonzero rows of $A$ by $\Gamma$, leading to a collection
$$
\Gcal=\{g_{A(\theta_j,\Gamma_k)}: 1\le j\le M_1,1\le k\le M_2\},
$$
where $A(\theta,\Gamma)$ is a $d\times N$ matrix whose $i$th row is zero if the $i$th entry of $\theta$ is zero, and the collection of the nonzero rows of $A$ is given by $\Gamma$. In what follows, for brevity, we shall write
$$
\Gcal=\{g_{A_k}: 1\le k\le M\},
$$
where $M=M_1M_2$ and
$$\Acal=\{A_k: 1\le k\le M\}$$
is the collection of $d\times N$ matrices of the form $A(\theta_j,\Gamma_k)$. By (c) and (b'),
$$
\log M\ge {1\over 4}s\log(d/s)+{1\over 8}Ns.
$$

Note that, for any two matrices $A,B\in \{-1,0,1\}^{d\times N}$ such that $s_A=s_B=:s$, we have
\begin{eqnarray*}
\|g_A-g_B\|_{L_2(\Pi)}^2&=&N^{-1}s^{-2/q}\int_{\Xcal^d}\left(\sum_{j=1}^d\sum_{k=1}^N (a_{jk}-b_{jk})\lambda^{1/2}_{j,N+k}\varphi_{j,N+k}(x_j)\right)^2d\Pi((x_1,\ldots,x_d)^\top)\\
&\ge&\eta_q^{-1}N^{-1}s^{-2/q}\sum_{j=1}^d\left\|\sum_{k=1}^N (a_{jk}-b_{jk})\lambda^{1/2}_{j,N+k}\varphi_{j,N+k}\right\|^2_{L_2(\Pi_j)}\\
&=&\eta_q^{-1}N^{-1}s^{-2/q}\sum_{j=1}^d\sum_{k=1}^N \lambda_{j,N+k}(a_{jk}-b_{jk})^2
\end{eqnarray*}
where the inequality follows from (\ref{eq:reassumption}). By (\ref{eq:eigassumption}), this can be further lower-bounded by
\begin{eqnarray*}
\|g_A-g_B\|_{L_2(\Pi)}^2&\ge&c_\lambda^{-1}\eta_q^{-1}N^{-1}s^{-2/q}\sum_{j=1}^d\sum_{k=1}^N \lambda_{N+k}(a_{jk}-b_{jk})^2\\
&\ge&c_\lambda^{-1}\eta_q^{-1}N^{-1}s^{-2/q}\lambda_{2N}\sum_{j=1}^d\sum_{k=1}^N (a_{jk}-b_{jk})^2\\
&=&c_\lambda^{-1}\eta_q^{-1}2^{-2\alpha}N^{-1-2\alpha}s^{-2/q}\|A-B\|_{\rm F}^2.
\end{eqnarray*}
By construction, for any $A\neq A'\in \Acal$,
$$
\|A-A'\|_{\rm F}^2\ge Ns/2,
$$
and hence,
$$
\|g_{A}-g_{A'}\|_{L_2(\Pi)}^2\ge c_\lambda^{-1}\eta_q^{-1}2^{-1-2\alpha}N^{-2\alpha}s^{1-2/q}.
$$

On the other hand, for any $A\in \Acal$,
\begin{eqnarray*}
\|g_A\|_{L_2(\Pi)}^2&=&N^{-1}s^{-2/q}\int_{\Xcal^d}\left(\sum_{j=1}^d\sum_{k=1}^N a_{jk}\lambda^{1/2}_{j,N+k}\varphi_{j,N+k}(x_j)\right)^2d\Pi((x_1,\ldots,x_d)^\top)\\
&\le&\eta_qN^{-1}s^{-2/q}\sum_{j=1}^d\left\|\sum_{k=1}^N a_{jk}\lambda^{1/2}_{j,N+k}\varphi_{j,N+k}\right\|^2_{L_2(\Pi_j)}\\
&=&\eta_qN^{-1}s^{-2/q}\sum_{j=1}^d\sum_{k=1}^N \lambda_{j,N+k}a_{jk}^2\\
&\le&c_\lambda\eta_qN^{-1}s^{-2/q}\sum_{j=1}^d\sum_{k=1}^N \lambda_{N+k}a_{jk}^2\\
&\le&c_\lambda\eta_qN^{-1}s^{-2/q}\lambda_N\sum_{j=1}^d\sum_{k=1}^N a_{jk}^2\\
&=&c_\lambda\eta_qN^{-2\alpha}s^{1-2/q}.
\end{eqnarray*}

Following a standard argument, the lower bound can be reduced to the error probability in a multi-way hypothesis test. See, e.g., Tsybakov (2009). More specifically, let $\Theta$ be a random variable uniformly distributed on $\{1,\ldots, M\}$. Then it can be deduced that
$$
\inf_{\tilde{f}}\sup_{f\in \Bcal_1(\ell_q(\Hcal_d))}\PP\left\{\|\tilde{f}-f\|_{L_2(\Pi)}^2\ge {1\over 4}\min_{A\neq A'\in\Acal}\|g_A-g_{A'}\|_{L_2(\Pi)}^2\right\}\ge \inf_{\hat{\Theta}}\PP\{\hat{\Theta}\neq \Theta\},
$$
where the infimum on the righthand side is taken over all decision rules that are measurable functions of the data. By Fano's Lemma, we get
\begin{equation}
\label{eq:fano}
\PP\left\{\hat{\Theta}\neq \Theta|X_1,\ldots,X_n\right\}\ge 1-{1\over \log M}\left[\II_{X_1,\ldots,X_n}(Y_1,\ldots,Y_n;\Theta)+\log 2\right],
\end{equation}
where $\II_{X_1,\ldots,X_n}(Y_1,\ldots,Y_n; \Theta)$ is the mutual information between $\Theta$ and $Y_1,\ldots,Y_n$ with $X_1,\ldots,X_n$ being held fixed. It is not hard to derive
\begin{eqnarray*}
\EE_{X_1,\ldots,X_n}\left[\II_{X_1,\ldots,X_n}(Y_1,\ldots,Y_n;\Theta)\right]&\le&\left(\begin{array}{c}M\\ 2\end{array}\right)^{-1}\sum_{A\neq A'\in \Acal}\EE_{X_1,\ldots,X_n}\Kcal(\bfP_{g_A}||\bfP_{g_{A'}})\\
&\le&{n\over 2}\left(\begin{array}{c}M\\ 2\end{array}\right)^{-1}\sum_{A\neq A'\in \Acal}\EE_{X_1,\ldots,X_n}\|g_{A}-g_{A'}\|_{L_2(\Pi_n)}^2,
\end{eqnarray*}
where $\Kcal(\cdot||\cdot)$ denote the Kullback-Leibler distance, $\bfP_g$ stands for conditional distribution of $\{Y_i: 1\le i\le n\}$ given $\{X_i: 1\le i\le n\}$ and the true regression function in (\ref{eq:regress}) is given by $f=g$, and for any $g: \Xcal^d\to \RR$,
$$
\|g\|_{L_2(\Pi_n)}^2={1\over n}\sum_{i=1}^n [g(X_i)]^2.
$$
Thus,
\begin{eqnarray*}
\EE_{X_1,\ldots,X_n}\left[\II_{X_1,\ldots,X_n}(Y_1,\ldots,Y_n;\Theta)\right]&\le&{n\over 2}\left(\begin{array}{c}M\\ 2\end{array}\right)^{-1}\sum_{A\neq A'\in \Acal}\|g_{A}-g_{A'}\|_{L_2(\Pi)}^2\\
&\le&{n\over 2}\max_{A\neq A'\in \Acal}\|g_{A}-g_{A'}\|_{L_2(\Pi)}^2\\
&\le&{2n}\max_{A\in \Acal}\|g_{A}\|_{L_2(\Pi)}^2\\
&\le&2c_\lambda\eta_qnN^{-2\alpha}s^{1-2/q}.
\end{eqnarray*}

Now, from (\ref{eq:fano}), we get
\begin{eqnarray*}
&&\inf_{\tilde{f}}\sup_{f\in \Bcal_1(\ell_q(\Hcal_d))}\PP\left\{\|\tilde{f}-f\|_2^2\ge c_\lambda^{-1}\eta_q^{-1}2^{-2-2\alpha}N^{-2\alpha}s^{1-2/q}\right\}\\
&\ge& \inf_{\hat{\Theta}}\PP\{\hat{\Theta}\neq \Theta\}\\
&\ge& 1-{\EE_{X_1,\ldots,X_n}\left[\II_{X_1,\ldots,X_n}(Y_1,\ldots,Y_n;\Theta)\right]+\log 2\over \log M}\\
&\ge&1-{2c_\lambda \eta_qnN^{-2\alpha}s^{1-2/q}+\log 2\over {{1\over 4}s\log (d/s)+{1\over 8}Ns}}.
\end{eqnarray*}

Taking $N=1$ and
$$
s=C_1 \left(n\over \log d\right)^{q/2}
$$
for a sufficiently small constant $C_1>0$ yields
\begin{equation}
\label{eq:lower1}
\inf_{\tilde{f}}\sup_{f\in \Bcal_1(\ell_q(\Hcal_d))}\PP\left\{\|\tilde{f}-f\|_2^2\ge C_2 \left(\log d\over n\right)^{1-q/2}\right\}\ge 3/4,
\end{equation}
for some constant $C_2>0$ depending on $\alpha$, $\eta_q$ and $c_\lambda$ only. On the other hand, if $\alpha\le 1/q-1/2$, taking
$$
s=1,\qquad {\rm and}\qquad N=C_1n^{1\over 2\alpha+1}
$$
for a sufficiently small constant $C_1>0$ yields
\begin{equation}
\label{eq:lower2}
\inf_{\tilde{f}}\sup_{f\in \Bcal_1(\ell_q(\Hcal_d))}\PP\left\{\|\tilde{f}-f\|_2^2\ge C_2 n^{-{2\alpha\over 2\alpha+1}}\right\}\ge 3/4.
\end{equation}
Combining (\ref{eq:lower1}) and (\ref{eq:lower2}), we have
$$
\inf_{\tilde{f}}\sup_{f\in \Bcal_1(\ell_q(\Hcal_d))}\PP\left\{\|\tilde{f}-f\|_2^2\ge C_2 \left[\left(\log d\over n\right)^{1-q/2}+n^{-{2\alpha\over 2\alpha+1}}\right]\right\}\ge 3/4,
$$
which completes the proof.

\subsection{Upper bounds}

We now prove the upper bounds given in Theorem \ref{th:upper}. By definition,
$$
{1\over n}\sum_{i=1}^n \left[Y_i-\hat{f}(X_i)\right]^2\le {1\over n}\sum_{i=1}^n \left[Y_i-f(X_i)\right]^2,
$$
which immediately implies that
\begin{equation}
\label{eq:basic}
{1\over n}\sum_{i=1}^n \left[\hat{f}(X_i)-f(X_i)\right]^2\le {2\over n}\sum_{i=1}^n\epsilon_i\left[\hat{f}(X_i)-f(X_i)\right].
\end{equation}
Write $\Delta_j=\hat{f}_j-f_j$ and $\Delta=\hat{f}-f$. It is clear that $\Delta=\sum_{j=1}^d \Delta_j$.

Our main strategy is to derive upper and lower bounds for the right and left hand side of (\ref{eq:basic}) respectively, and then put them together to derive (\ref{eq:bdlse}).

\paragraph{Step 1. Bounding the righthand side of (\ref{eq:basic}).}

Observe that
$$
\left|{1\over n}\sum_{i=1}^n \epsilon_i \Delta_j(x_{ij})\right|\le \|\Delta_j\|_{\Hcal_1}\hat{Z}_{jn}\left({\|\Delta_j\|_{L_2(\Pi_{jn})}\over\|\Delta_j\|_{\Hcal_1}}\right),
$$
where $\hat{Z}_{jn}$ is defined by (\ref{eq:gausscomp}). By Lemma \ref{le:emp}, this can be further bounded by
$$
C_1n^{-1/2}\left(\|\Delta_j\|_{L_2(\Pi_{jn})}^{1-{1\over 2\alpha}}\|\Delta_j\|_{\Hcal_1}^{1\over 2\alpha}+\|\Delta_j\|_{L_2(\Pi_{jn})}\sqrt{(\beta+1)\log d}+e^{-d}\|\Delta_j\|_{\Hcal_1}\right)
$$
for some constant $C_1>0$, with probability at least $1-d^{-(\beta+1)}$. By union bound, with probability $1-d^{-\beta}$,
\begin{eqnarray}
\nonumber
{2\over n}\sum_{i=1}^n\epsilon_i\left[\hat{f}(X_i)-f(X_i)\right]&\le& 2\sum_{j=1}^d\left|{1\over n}\sum_{i=1}^n \epsilon_i \Delta_j(x_{ij})\right|\\
\nonumber
&\le&2C_1n^{-1/2}\sum_{j=1}^d\|\Delta_j\|_{L_2(\Pi_{jn})}^{1-{1\over 2\alpha}}\|\Delta_j\|_{\Hcal_1}^{1\over 2\alpha}\\
\nonumber
&&\qquad +2C_1n^{-1/2}\sqrt{(\beta+1)\log d}\sum_{j=1}^d\|\Delta_j\|_{L_2(\Pi_{jn})}\\
&&\qquad +2C_1n^{-1/2}e^{-d}\sum_{j=1}^d\|\Delta_j\|_{\Hcal_1}.
\label{eq:sepbds}
\end{eqnarray}
We denote by $\Ecal_1$ the event that the above inequality holds. We now bound the three terms on the rightmost side separately.

We first derive a bound for
$$
n^{-1/2}\sum_{j=1}^d\|\Delta_j\|_{\Hcal_1}^{1\over 2\alpha}\|\Delta_j\|_{L_2(\Pi_{jn})}^{1-{1\over 2\alpha}}.
$$
We treat the cases of $2/(2\alpha+1)\ge q$ and $2/(2\alpha+1)< q$ separately.

\begin{enumerate}
\item[Case 1:] $2/(2\alpha+1)\ge q$. By Young's inequality, for a constant $\zeta>1$ whose value will be specified later,
$$
n^{-1/2}\|\Delta_j\|_{\Hcal_1}^{1\over 2\alpha}\|\Delta_j\|_{L_2(\Pi_{jn})}^{1-{1\over 2\alpha}}\le \zeta^{-{4\alpha\over 2\alpha-1}}\|\Delta_j\|_{L_2(\Pi_{jn})}^2+\zeta^{4\alpha\over 2\alpha+1}n^{-{2\alpha\over 2\alpha+1}}\|\Delta_j\|_{\Hcal_1}^{2\over 2\alpha+1}.
$$
Note that for any $q\le q'\le 2$,
\begin{eqnarray*}
\sum_{j=1}^d\|\Delta_j\|_{\Hcal_1}^{q'}&\le& 2\left(\sum_{j=1}^d\|\hat{f}_j\|_{\Hcal_1}^{q'}+\sum_{j=1}^d\|f_j\|_{\Hcal_1}^{q'}\right)\\
&\le&2\left(\sum_{j=1}^d\|\hat{f}_j\|_{\Hcal_1}^q+\sum_{j=1}^d\|f_j\|_{\Hcal_1}^q\right)\\
&\le&4.
\end{eqnarray*}
In particular, we get
$$
\sum_{j=1}^d \|\Delta_j\|_{\Hcal_1}^{2\over 2\alpha+1}\le 4,
$$
Hence,
\begin{equation}
\label{eq:bde1}
\sum_{j=1}^d n^{-1/2}\|\Delta_j\|_{\Hcal_1}^{1\over 2\alpha}\|\Delta_j\|_{L_2(\Pi_{jn})}^{1-{1\over 2\alpha}}\le \zeta^{-{4\alpha\over 2\alpha-1}}\|\Delta_j\|_{L_2(\Pi_{jn})}^2+4\zeta^{4\alpha\over 2\alpha+1}n^{-{2\alpha\over 2\alpha+1}}.
\end{equation}

\item[Case 2:] $2/(2\alpha+1)<q$.Write
\begin{eqnarray*}
&&n^{-1/2}\sum_{j=1}^d\|\Delta_j\|_{\Hcal_1}^{1\over 2\alpha}\|\Delta_j\|_{L_2(\Pi_{jn})}^{1-{1\over 2\alpha}}\\
&=&n^{-1/2}\sum_{j: \|\Delta_j\|_{\Hcal_1}>n^{-1/2}}\|\Delta_j\|_{\Hcal_1}^{1\over 2\alpha}\|\Delta_j\|_{L_2(\Pi_{jn})}^{1-{1\over 2\alpha}}\\
&&\qquad\qquad +n^{-1/2}\sum_{j: \|\Delta_j\|_{\Hcal_1}\le n^{-1/2}}\|\Delta_j\|_{\Hcal_1}^{1\over 2\alpha}\|\Delta_j\|_{L_2(\Pi_{jn})}^{1-{1\over 2\alpha}}.
\end{eqnarray*}
For the first term on the right hand side, by a similar argument as before, we have
\begin{eqnarray*}
&&n^{-1/2}\sum_{j: \|\Delta_j\|_{\Hcal_1}>n^{-1/2}}\|\Delta_j\|_{\Hcal_1}^{1\over 2\alpha}\|\Delta_j\|_{L_2(\Pi_{jn})}^{1-{1\over 2\alpha}}\\
&\le& \zeta^{-{4\alpha\over 2\alpha-1}}\sum_{j: \|\Delta_j\|_{\Hcal_1}>n^{-1/2}}\|\Delta_j\|_{L_2(\Pi_{jn})}^2+\zeta^{{4\alpha\over 2\alpha+1}}n^{-{2\alpha\over 2\alpha+1}}\sum_{j: \|\Delta_j\|_{\Hcal_1}>n^{-1/2}}\|\Delta_j\|_{\Hcal_1}^{2\over 2\alpha+1}\\
&\le& \zeta^{-{4\alpha\over 2\alpha-1}}\sum_{j: \|\Delta_j\|_{\Hcal_1}>n^{-1/2}}\|\Delta_j\|_{L_2(\Pi_{jn})}^2+\zeta^{{4\alpha\over 2\alpha+1}}n^{-(1-{q\over 2})}\sum_{j: \|\Delta_j\|_{\Hcal_1}>n^{-1/2}}\|\Delta_j\|_{\Hcal_1}^{q}\\
&\le&\zeta^{-{4\alpha\over 2\alpha-1}}\sum_{j: \|\Delta_j\|_{\Hcal_1}>n^{-1/2}}\|\Delta_j\|_{L_2(\Pi_{jn})}^2+4\zeta^{{4\alpha\over 2\alpha+1}}n^{-(1-{q\over 2})},
\end{eqnarray*}
where in the last inequality, we used the fact that
$$
\sum_{j: \|\Delta_j\|_{\Hcal_1}>n^{-1/2}}\|\Delta_j\|_{\Hcal_1}^{q}\le \sum_{j=1}^d\|\Delta_j\|_{\Hcal_1}^{q}\le 2\sum_{j=1}^d\left(\|\hat{f}_j\|_{\Hcal_1}^{q}+\|{f}_j\|_{\Hcal_1}^{q}\right)\le 4.
$$
On the other hand, because
$$\|\Delta_j\|_{L_2(\Pi_{jn})}\le \|\Delta_j\|_{L_\infty}\le \|\Delta_j\|_{\Hcal_1},$$
we get
\begin{eqnarray*}
n^{-1/2}\sum_{j: \|\Delta_j\|_{\Hcal_1}\le n^{-1/2}}\|\Delta_j\|_{\Hcal_1}^{1\over 2\alpha}\|\Delta_j\|_{L_2(\Pi_{jn})}^{1-{1\over 2\alpha}}&\le&n^{-1/2}\sum_{j: \|\Delta_j\|_{\Hcal_1}\le n^{-1/2}}\|\Delta_j\|_{\Hcal_1}\\
&\le&n^{-(1-q/2)}\sum_{j: \|\Delta_j\|_{\Hcal_1}\le n^{-1/2}}\|\Delta_j\|_{\Hcal_1}^q\\
&\le&n^{-(1-q/2)}\sum_{j=1}^d\|\Delta_j\|_{\Hcal_1}^q\\
&\le&4n^{-(1-q/2)}.
\end{eqnarray*}
Thus,
\begin{equation}
\label{eq:bde2}
n^{-1/2}\sum_{j=1}^d\|\Delta_j\|_{\Hcal_1}^{1\over 2\alpha}\|\Delta_j\|_{L_2(\Pi_{jn})}^{1-{1\over 2\alpha}}\le \zeta^{-{4\alpha\over 2\alpha-1}}\sum_{j=1}^d\|\Delta_j\|_{L_2(\Pi_{jn})}^2+8\zeta^{{4\alpha\over 2\alpha+1}}n^{-(1-{q\over 2})}.
\end{equation}
\end{enumerate}
Combing (\ref{eq:bde1}) and (\ref{eq:bde2}), we get
\begin{equation}
\label{eq:bde3}
n^{-1/2}\sum_{j=1}^d\|\Delta_j\|_{\Hcal_1}^{1\over 2\alpha}\|\Delta_j\|_{L_2(\Pi_{jn})}^{1-{1\over 2\alpha}}\le \zeta^{-{4\alpha\over 2\alpha-1}}\sum_{j=1}^d\|\Delta_j\|_{L_2(\Pi_{jn})}^2+8\zeta^{{4\alpha\over 2\alpha+1}}n^{-(1-\max\{{q\over 2},{1\over 2\alpha+1}\})}.
\end{equation}

By Theorem 4 of Koltchinskii and Yuan (2010), there exists a numerical constant $C_2>1$ such that with probability at least $1-d^{-\beta}$ for all $h\in \Hcal_1$, and $j=1,\ldots, d$,
\begin{equation}
\label{eq:nbd1}
\|h\|_{L_2(\Pi_j)}\le C_2\left[\|h\|_{L_2(\Pi_{jn})}+\left(n^{-{\alpha\over 2\alpha+1}}+\sqrt{(\beta+1)\log d\over n}\right)\|h\|_{\Hcal_1}\right],
\end{equation}
and
\begin{equation}
\label{eq:nbd2}
\|h\|_{L_2(\Pi_{jn})}\le C_2\left[\|h\|_{L_2(\Pi_j)}+\left(n^{-{\alpha\over 2\alpha+1}}+\sqrt{(\beta+1)\log d\over n}\right)\|h\|_{\Hcal_1}\right].
\end{equation}
Denote by $\Ecal_2$ the event that both (\ref{eq:nbd1}) and (\ref{eq:nbd2}) hold. Under $\Ecal_2$,
\begin{eqnarray*}
\sum_{j=1}^d\|\Delta_j\|_{L_2(\Pi_{jn})}^2&\le& 2C_2^2\sum_{j=1}^d\left[\|\Delta_j\|^2_{L_2(\Pi_j)}+\left(n^{-{2\alpha\over 2\alpha+1}}+{(\beta+1)\log d\over n}\right)\|\Delta_j\|_{\Hcal_1}^2\right]\\
&\le&2C_2^2\sum_{j=1}^d\|\Delta_j\|^2_{L_2(\Pi_j)}+8C_2^2\left(n^{-{2\alpha\over 2\alpha+1}}+{(\beta+1)\log d\over n}\right),
\end{eqnarray*}
where the second inequality follows from the fact that
$$
\sum_{j=1}^d\|\Delta_j\|_{\Hcal_1}^2\le 4.
$$
By (\ref{eq:reassumption}), this implies that
$$
\sum_{j=1}^d\|\Delta_j\|_{L_2(\Pi_{jn})}^2\le 2C_2^2\eta_q\|\Delta\|_{L_2(\Pi)}^2+8C_2^2\left(n^{-{2\alpha\over 2\alpha+1}}+{(\beta+1)\log d\over n}\right).
$$
Together with (\ref{eq:bde3}), we get
\begin{eqnarray}
\nonumber
n^{-1/2}\sum_{j=1}^d\|\Delta_j\|_{\Hcal_1}^{1\over 2\alpha}\|\Delta_j\|_{L_2(\Pi_{jn})}^{1-{1\over 2\alpha}}&\le& 2C_2^2\eta_q\zeta^{-{4\alpha\over 2\alpha-1}}\|\Delta\|_{L_2(\Pi)}^2\\
\nonumber
&&\qquad +8C_2^2\zeta^{-{4\alpha\over 2\alpha-1}}\left(n^{-{2\alpha\over 2\alpha+1}}+{(\beta+1)\log d\over n}\right)\\
&&\qquad +8\zeta^{{4\alpha\over 2\alpha+1}}n^{-(1-\max\{{q\over 2},{1\over 2\alpha+1}\})}.
\label{eq:bde4}
\end{eqnarray}

The second term on the rightmost hand side of (\ref{eq:sepbds}) can also be bounded under event $\Ecal_2$. By (\ref{eq:nbd2}),
\begin{eqnarray}
\nonumber
\sum_{j=1}^d\|\Delta_j\|_{L_2(\Pi_{jn})}&\le& C_2\sum_{j=1}^d\|\Delta_j\|_{L_2(\Pi_j)}+C_2\left(n^{-{\alpha\over 2\alpha+1}}+\sqrt{(\beta+1)\log d\over n}\right)\sum_{j=1}^d\|\Delta_j\|_{\Hcal_1}\\
&\le&C_2\sum_{j=1}^d\|\Delta_j\|_{L_2(\Pi_j)}+4C_2\left(n^{-{\alpha\over 2\alpha+1}}+\sqrt{(\beta+1)\log d\over n}\right),
\label{eq:bde5}
\end{eqnarray}
where in the second inequality we used the fact that
$$
\sum_{j=1}^d\|\Delta_j\|_{\Hcal_1}\le \sum_{j=1}^d\|\Delta_j\|_{\Hcal_1}^q\le 4.
$$
Write
$$
\sum_{j=1}^d\|\Delta_j\|_{L_2(\Pi_j)}\le \sum_{j:\|\Delta_j\|_{L_2(\Pi_j)}> \sqrt{\log d\over n}}\|\Delta_j\|_{L_2(\Pi_j)}+\sum_{j:\|\Delta_j\|_{L_2(\Pi_j)}\le \sqrt{\log d\over n}}\|\Delta_j\|_{L_2(\Pi_j)}.
$$
The first term can be bounded by Cachy-Schwartz inequality:
\begin{eqnarray*}
&&\sum_{j:\|\Delta_j\|_{L_2(\Pi_j)}> \sqrt{\log d\over n}}\|\Delta_j\|_{L_2(\Pi_j)}\\
&&\qquad \le \left({\rm card}\left\{j:\|\Delta_j\|_{L_2(\Pi_j)}> \sqrt{\log d\over n}\right\}\right)^{1/2}\left(\sum_{j:\|\Delta_j\|_{L_2(\Pi_j)}> \sqrt{\log d\over n}}\|\Delta_j\|_{L_2(\Pi_j)}^2\right)^{1/2}.
\end{eqnarray*}
Observe that
$$
{\rm card}\left\{j:\|\Delta_j\|_{L_2(\Pi_j)}> \sqrt{\log d\over n}\right\}\le \left({\log d\over n}\right)^{-q/2}\sum_{j=1}^d\|\Delta_j\|_{\Hcal_1}^q\le 4\left({\log d\over n}\right)^{-q/2}.
$$
Thus,
\begin{eqnarray*}
\sum_{j:\|\Delta_j\|_{L_2(\Pi_j)}> \sqrt{\log d\over n}}\|\Delta_j\|_{L_2(\Pi_j)}&\le& 4\left({\log d\over n}\right)^{-q/4}\left(\sum_{j=1}^d\|\Delta_j\|_{L_2(\Pi_j)}^2\right)^{1/2}\\
&\le&4\eta_q^{1/2}\left({\log d\over n}\right)^{-q/4}\|\Delta\|_{L_2(\Pi)}.
\end{eqnarray*}
Together with the fact that
\begin{eqnarray*}
\sum_{j:\|\Delta_j\|_{L_2(\Pi_j)}\le \sqrt{\log d\over n}}\|\Delta_j\|_{L_2(\Pi_j)}&\le& \sum_{j:\|\Delta_j\|_{L_2(\Pi_j)}\le \sqrt{\log d\over n}}\|\Delta_j\|_{L_2(\Pi_j)}^q\left({\log d\over n}\right)^{(1-q)/2}\\
&\le&\left({\log d\over n}\right)^{(1-q)/2}\sum_{j=1}^d\|\Delta_j\|_{L_2(\Pi_j)}^q\\
&\le&4\left({\log d\over n}\right)^{(1-q)/2},
\end{eqnarray*}
we get
\begin{equation}
\label{eq:bdl1}
\sum_{j=1}^d\|\Delta_j\|_{L_2(\Pi_j)}\le 4\eta_q^{1/2}\left({\log d\over n}\right)^{-q/4}\|\Delta\|_{L_2(\Pi)}+4\left({\log d\over n}\right)^{1-q\over 2}.
\end{equation}
In the light of (\ref{eq:bde5}), we have
\begin{eqnarray}
\nonumber
\sqrt{\log d\over n}\sum_{j=1}^d\|\Delta_j\|_{L_2(\Pi_{jn})}&\le& 4C_2\eta_q^{1/2} \left({\log d\over n}\right)^{1/2-q/4}\|\Delta\|_{L_2(\Pi)}\\
&&\qquad+4C_2n^{-{\alpha\over 2\alpha+1}}\sqrt{\log d\over n}+8C_2\sqrt{\beta+1}\left({\log d\over n}\right)^{1-{q\over 2}},
\label{eq:bde6}
\end{eqnarray}
where we used the fact that $\log d<n$ and $C_2>1$.

Combing (\ref{eq:sepbds}), (\ref{eq:bde4}), (\ref{eq:bde6}) and the fact that
$$
\sum_{j=1}^d\|\Delta_j\|_{\Hcal_1}\le 4,
$$
we get
\begin{eqnarray}
\nonumber
{2\over n}\sum_{i=1}^n\epsilon_i\left[\hat{f}(X_i)-f(X_i)\right]&\le&C_3\eta_q\zeta^{-{4\alpha\over 2\alpha-1}}\|\Delta\|_{L_2(\Pi)}^2\\
\nonumber
&&\qquad +C_3\zeta^{-{4\alpha\over 2\alpha-1}}\left(n^{-{2\alpha\over 2\alpha+1}}+{(\beta+1)\log d\over n}\right)\\
\nonumber
&&\qquad +C_3\zeta^{{4\alpha\over 2\alpha+1}}n^{-(1-\max\{{q\over 2},{1\over 2\alpha+1}\})}\\
\nonumber
&&\qquad +C_3\sqrt{\beta+1}\eta_q^{1/2} \left({\log d\over n}\right)^{1/2-q/4}\|\Delta\|_{L_2(\Pi)}\\
\nonumber
&&\qquad+C_3\sqrt{\beta+1}n^{-{\alpha\over 2\alpha+1}}\sqrt{\log d\over n}\\
\nonumber
&&\qquad +C_3\sqrt{\beta+1}\left({\log d\over n}\right)^{1-{q\over 2}}\\
&&\qquad +C_3n^{-1/2}e^{-d},
\label{eq:rhsbd}
\end{eqnarray}
for some constant $C_3>0$, under the event $\Ecal_1\cap \Ecal_2$.
\paragraph{Step 2. Bounding the lefthand side of (\ref{eq:basic}).}

To bound the lefthand side of (\ref{eq:basic}), first observe that
\begin{equation}
\label{eq:bdbasicleft}
\|\Delta\|_{L_2(\Pi)}^2-\|\Delta\|_{L_2(\Pi_n)}^2\le \sup_{\substack{g\in \Bcal_4(\ell_q(\Hcal_d))\\ \|g\|_{L_2(\Pi)}\le \|\Delta\|_{L_2(\Pi)}}}\left(\|g\|_{L_2(\Pi)}^2-\|g\|_{L_2(\Pi_n)}^2\right)
\end{equation}
Note that for any $g\in \Bcal_4(\ell_q(\Hcal_d))$,
$$
\|g\|_{L_\infty}^2\le \|g\|_{\ell_1(\Hcal_d)}^2\le \left(\|g\|_{\ell_q(\Hcal_d)}^q\right)^2\le 16,
$$
and
$$
\|g\|_{L_2(\Pi)}^4\le \|g\|_{L_\infty}^2\|g\|_{L_2(\Pi)}^2\le 16\|g\|_{L_2(\Pi)}^2.
$$
By Talagrand's concentration inequality, for any fixed $u\in [0,1]$,
\begin{eqnarray*}
&&\sup_{\substack{g\in \Bcal_4(\ell_q(\Hcal_d))\\ \|g\|_{L_2(\Pi)}\le u}}\left(\|g\|_{L_2(\Pi)}^2-\|g\|_{L_2(\Pi_n)}^2\right)\\
&\le& 2\left(\EE \sup_{\substack{g\in \Bcal_4(\ell_q(\Hcal_d))\\ \|g\|_{L_2(\Pi)}\le u}}\left(\|g\|_{L_2(\Pi)}^2-\|g\|_{L_2(\Pi_n)}^2\right)+4u\sqrt{t\over n}+{16t\over n}\right),
\end{eqnarray*}
with probability at least $1-e^{-t}$. By symmetrization inequality,
$$
\EE \sup_{\substack{g\in \Bcal_4(\ell_q(\Hcal_d))\\ \|g\|_{L_2(\Pi)}\le u}}\left(\|g\|_{L_2(\Pi)}^2-\|g\|_{L_2(\Pi_n)}^2\right)\le 2\EE \sup_{\substack{g\in \Bcal_4(\ell_q(\Hcal_d))\\ \|g\|_{L_2(\Pi)}\le u}}\left({1\over n}\sum_{i=1}^n \sigma_i g^2(X_i)\right).
$$
Note that $g^2$ is 8-Lipschitz function on $\Bcal_4(\ell_q(\Hcal_d))$. By contraction inequality,
$$
\EE \sup_{\substack{g\in \Bcal_4(\ell_q(\Hcal_d))\\ \|g\|_{L_2(\Pi)}\le u}}\left({1\over n}\sum_{i=1}^n \sigma_i g^2(X_i)\right)\le 8\EE \sup_{\substack{g\in \Bcal_4(\ell_q(\Hcal_d))\\ \|g\|_{L_2(\Pi)}\le u}}\left({1\over n}\sum_{i=1}^n \sigma_i g(X_i)\right).
$$
Again by Talagrand's concentration inequality, there exists a numerical constant $C_4>0$ such that with probability at least $1-e^{-t}$,
\begin{eqnarray*}
&&\EE \sup_{\substack{g\in \Bcal_4(\ell_q(\Hcal_d))\\ \|g\|_{L_2(\Pi)}\le u}}\left({1\over n}\sum_{i=1}^n \sigma_i g(X_i)\right)\\
&\le& C_4\left(\sup_{\substack{g\in \Bcal_4(\ell_q(\Hcal_d))\\ \|g\|_{L_2(\Pi)}\le u}}\left({1\over n}\sum_{i=1}^n \sigma_i g(X_i)\right)+u\sqrt{t\over n}+{t\over n}\right)\\
&\le& C_4\left(\sup_{\substack{\sum_{j=1}^d\|g_j\|_{\Hcal_1}^q\le 4\\ \left\|\sum_{j=1}^dg_j\right\|_{L_2(\Pi)}\le u}}\sum_{j=1}^d\left({1\over n}\sum_{i=1}^n \sigma_i g_j(x_{ij})\right)+u\sqrt{t\over n}+{t\over n}\right).
\end{eqnarray*}
In other words,
\begin{eqnarray}
\nonumber
&&\sup_{\substack{g\in \Bcal_4(\ell_q(\Hcal_d))\\ \|g\|_{L_2(\Pi)}\le u}}\left(\|g\|_{L_2(\Pi)}^2-\|g\|_{L_2(\Pi_n)}^2\right)\\
&\le& 16C_4\left(\sup_{\substack{\sum_{j=1}^d\|g_j\|_{\Hcal_1}^q\le 4\\ \left\|\sum_{j=1}^dg_j\right\|_{L_2(\Pi)}\le u}}\sum_{j=1}^d\left({1\over n}\sum_{i=1}^n \sigma_i g_j(x_{ij})\right)+u\sqrt{t\over n}+{t\over n}\right),
\label{eq:bdsq1}
\end{eqnarray}
with probability at least $1-2e^{-t}$.

Note that
$$
{1\over n}\sum_{i=1}^n \sigma_i g_j(x_{ij})\le \|g_j\|_{\Hcal_1} \sup_{\substack{\|h\|_{\Hcal_1}=1\\ \|h\|_{L_2(\Pi_j)}\le \|g_j\|_{L_2(\Pi_j)}/\|g_j\|_{\Hcal_1}}}\left({1\over n}\sum_{i=1}^n \sigma_i h(x_{ij})\right)
$$
By Lemma \ref{le:emp} and union bound, there exists a constant $C_5>0$ such that
$$
\sup_{\substack{\|h\|_{\Hcal_1}=1\\ \|h\|_{L_2(\Pi_j)}\le u}}\left({1\over n}\sum_{i=1}^n \sigma_i h(x_{ij})\right)\le C_5n^{-1/2}\left(u^{1-{1\over 2\alpha}}+u\sqrt{(\beta+1)\log d}+e^{-d}\right),
$$
uniformly over $u\in [0,1]$ and $j=1,\ldots,d$ with probability at least $1-d^{-\beta}$. Denote this event by $\Ecal_3$, and we shall now proceed conditional on $\Ecal_3$.

It is not hard to see that, under $\Ecal_3$,
\begin{eqnarray}
\nonumber
&&\sum_{j=1}^d\left({1\over n}\sum_{i=1}^n\sigma_i g_j(x_{ij})\right)\\
&\le& C_5n^{-1/2}\sum_{j=1}^d\left(\|g_j\|^{1\over 2\alpha}_{\Hcal_1}\|g_j\|^{1-{1\over 2\alpha}}_{L_2(\Pi_j)}+\|g_j\|_{L_2(\Pi_j)}\sqrt{(\beta+1)\log d}+e^{-d}\|g_j\|_{\Hcal_1}\right).
\label{eq:bdrad1}
\end{eqnarray}
Following the same argument as that for (\ref{eq:bde3}), it can derived
\begin{eqnarray}
\nonumber
&&n^{-1/2}\sup_{\substack{\sum_{j=1}^d\|g_j\|_{\Hcal_1}^q\le 4\\ \left\|\sum_{j=1}^dg_j\right\|_{L_2(\Pi)}\le u}} \sum_{j=1}^d\|g_j\|_{\Hcal_1}^{1\over 2\alpha}\|g_j\|_{L_2(\Pi_j)}^{1-{1\over 2\alpha}}\\
\nonumber
&\le& \zeta^{-{4\alpha\over 2\alpha-1}}\sup_{\substack{\sum_{j=1}^d\|g_j\|_{\Hcal_1}^q\le 4\\ \left\|\sum_{j=1}^dg_j\right\|_{L_2(\Pi)}\le u}}\sum_{j=1}^d\|g_j\|_{L_2(\Pi_j)}^2+8\zeta^{{4\alpha\over 2\alpha+1}}n^{-(1-\max\{{q\over 2},{1\over 2\alpha+1}\})}\\
&\le&\zeta^{-{4\alpha\over 2\alpha-1}}\eta_qu^2+8\zeta^{{4\alpha\over 2\alpha+1}}n^{-(1-\max\{{q\over 2},{1\over 2\alpha+1}\})}.
\label{eq:bdrad2}
\end{eqnarray}
Similar to (\ref{eq:bdl1}), it can also be shown that for any $g_1,\ldots, g_d$ such that
$$
\sum_{j=1}^d\|g_j\|_{\Hcal_1}^q\le 4\qquad {\rm and}\qquad \sum_{j=1}^d\|g_j\|_{L_2(\Pi_j)}\le u,
$$
we have
\begin{equation}
\label{eq:bdrad3}
\sum_{j=1}^d\|g_j\|_{L_2(\Pi_j)}\le 4\eta_q^{1/2}\left({\log d\over n}\right)^{-q/4}u+4\left({\log d\over n}\right)^{1-q\over 2}.
\end{equation}
Combining (\ref{eq:bdrad1}), (\ref{eq:bdrad2}) and (\ref{eq:bdrad3}), we have
\begin{eqnarray*}
&&\sup_{\substack{\sum_{j=1}^d\|g_j\|_{\Hcal_1}^q\le 4\\ \left\|\sum_{j=1}^dg_j\right\|_{L_2(\Pi)}\le u}}\sum_{j=1}^d\left({1\over n}\sum_{i=1}^n \sigma_i g_j(x_{ij})\right)\\
&\le& C_5\zeta^{-{4\alpha\over 2\alpha-1}}\eta_qu^2+8C_5\zeta^{{4\alpha\over 2\alpha+1}}n^{-(1-\max\{{q\over 2},{1\over 2\alpha+1}\})}\\
&&\quad +4C_5\sqrt{(\beta+1)\log d\over n}\left(\eta_q^{1/2}\left({\log d\over n}\right)^{-q/4}u+\left({\log d\over n}\right)^{1-q\over 2}\right)\\
&&\quad +C_5n^{-1/2}e^{-d}.
\end{eqnarray*}
Together with (\ref{eq:bdsq1}), conditional on $\Ecal_3$,
\begin{eqnarray*}
&&\sup_{\substack{g\in \Bcal_4(\ell_q(\Hcal_d))\\ \|g\|_{L_2(\Pi)}\le u}}\left(\|g\|_{L_2(\Pi)}^2-\|g\|_{L_2(\Pi_n)}^2\right)\\
&\le& C_6\zeta^{-{4\alpha\over 2\alpha-1}}\eta_qu^2+C_6\zeta^{{4\alpha\over 2\alpha+1}}n^{-(1-\max\{{q\over 2},{1\over 2\alpha+1}\})}\\
&&\qquad +C_6\sqrt{(\beta+1)\log d\over n}\left(\eta_q^{1/2}\left({\log d\over n}\right)^{-q/4}u+\left({\log d\over n}\right)^{1-q\over 2}\right)\\
&&\qquad +C_6n^{-1/2}e^{-d}+C_6\left(u\sqrt{t\over n}+{t\over n}\right).
\end{eqnarray*}
holds for some constant $C_6>0$, with probability at least $1-2e^{-t}$. Using a peeling argument similar to that for Lemma \ref{le:emp1}, we can make this bound uniformly over $u\in [0,1]$. More specifically, it can be shown that there exist constants $C_7>0$ such that, conditional on $\Ecal_3$,
\begin{eqnarray}
\nonumber
&&\sup_{\substack{g\in \Bcal_4(\ell_q(\Hcal_d))\\ \|g\|_{L_2(\Pi)}\le u}}\left(\|g\|_{L_2(\Pi)}^2-\|g\|_{L_2(\Pi_n)}^2\right)\\
\nonumber
&\le& C_6\zeta^{-{4\alpha\over 2\alpha-1}}\eta_qu^2+C_6\zeta^{{4\alpha\over 2\alpha+1}}n^{-(1-\max\{{q\over 2},{1\over 2\alpha+1}\})}\\
\nonumber
&&\qquad +C_6\sqrt{(\beta+1)\log d\over n}\left(\eta_q^{1/2}\left({\log d\over n}\right)^{-q/4}u+\left({\log d\over n}\right)^{1-q\over 2}\right)\\
\nonumber
&&\qquad +C_6n^{-1/2}e^{-d}\\
&&\qquad +C_7\left(u\sqrt{(\beta+1)\log d\over n}+{(\beta+1)\log d\over n}\right),
\label{eq:bdleft}
\end{eqnarray}
uniformly over all $u\in[0,1]$ with probability at least $1-d^{-\beta}$. Denote by $\Ecal_4$ the event that inequality (\ref{eq:bdleft}) holds. Then
$$
\PP\{\Ecal_4\}\ge \PP\{\Ecal_4|\Ecal_3\}\PP(\Ecal_3)\ge (1-d^{-\beta})^2\ge 1-2d^{-\beta}.
$$

Together with (\ref{eq:bdbasicleft}), we get, under event $\Ecal_4$, 
\begin{eqnarray}
\nonumber
\|\Delta\|_{L_2(\Pi)}^2&\le& \|\Delta\|_{L_2(\Pi_n)}^2+C_8\zeta^{-{4\alpha\over 2\alpha-1}}\eta_qu^2+C_8\zeta^{{4\alpha\over 2\alpha+1}}n^{-(1-\max\{{q\over 2},{1\over 2\alpha+1}\})}\\
\nonumber
&&\qquad +C_8\sqrt{(\beta+1)\log d\over n}\left(\eta_q^{1/2}\left({\log d\over n}\right)^{-q/4}\|\Delta\|_{L_2(\Pi)}+\left({\log d\over n}\right)^{1-q\over 2}\right)\\
\nonumber
&&\qquad +C_8n^{-1/2}e^{-d}\\
&&\qquad +C_8\left(u\sqrt{(\beta+1)\log d\over n}+{(\beta+1)\log d\over n}\right),
\label{eq:bdleftfinal}
\end{eqnarray}
for some constant $C_8>0$.
\paragraph{Step 3. Putting it together.} Combining (\ref{eq:rhsbd}) and (\ref{eq:bdleftfinal}), we get
\begin{eqnarray*}
\|\Delta\|_{L_2(\Pi)}^2&\le&C_9\eta_q\zeta^{-{4\alpha\over 2\alpha-1}}\|\Delta\|_{L_2(\Pi)}^2\\
&&\qquad +C_9\zeta^{-{4\alpha\over 2\alpha-1}}{(\beta+1)\log d\over n}\\
&&\qquad +C_9\zeta^{{4\alpha\over 2\alpha+1}}n^{-(1-\max\{{q\over 2},{1\over 2\alpha+1}\})}\\
&&\qquad +C_9\sqrt{\beta+1}\left({\log d\over n}\right)^{1/2-q/4}\|\Delta\|_{L_2(\Pi)}\\
&&\qquad+C_9\sqrt{\beta+1}n^{-{\alpha\over 2\alpha+1}}\sqrt{\log d\over n}\\
&&\qquad +C_9(\beta+1)\left({\log d\over n}\right)^{1-{q\over 2}}\\
&&\qquad +C_9n^{-1/2}e^{-d},
\end{eqnarray*}
for some constant $C_9>0$, under the event $\Ecal_1\cap \Ecal_2\cap \Ecal_4$.

Take $\zeta$ large enough so that 
$$
C_9\eta_q\zeta^{-{4\alpha\over 2\alpha-1}}\le 1/2.
$$
Then
\begin{eqnarray*}
\|\Delta\|_{L_2(\Pi)}^2&\le& 2C_9\zeta^{-{4\alpha\over 2\alpha-1}}{(\beta+1)\log d\over n}\\
&&\qquad +2C_9\zeta^{{4\alpha\over 2\alpha+1}}n^{-(1-\max\{{q\over 2},{1\over 2\alpha+1}\})}\\
&&\qquad +2C_9\sqrt{\beta+1}\left({\log d\over n}\right)^{1/2-q/4}\|\Delta\|_{L_2(\Pi)}\\
&&\qquad+2C_9\sqrt{\beta+1}n^{-{\alpha\over 2\alpha+1}}\sqrt{\log d\over n}\\
&&\qquad +2C_9\sqrt{\beta+1}\left({\log d\over n}\right)^{1-{q\over 2}}\\
&&\qquad +2C_9n^{-1/2}e^{-d}.
\end{eqnarray*}
Therefore, there exists a constant $C_{10}>0$ such that, under the event $\Ecal_1\cap \Ecal_2\cap \Ecal_4$,
$$
\|\Delta\|_{L_2(\Pi)}^2\le C_{10}(\beta+1)\left(n^{-{2\alpha\over 2\alpha+1}}+\left({\log d\over n}\right)^{1-{q\over 2}}+\left({\log d\over n}\right)^{1/2-q/4}\|\Delta\|_{L_2(\Pi)}\right),
$$
which implies (\ref{eq:bdlse}). Statement (\ref{eq:bdlse}) now follows from the fact that
$$
\PP\{\Ecal_1\cap\Ecal_2\cap\Ecal_4\}\ge 1-\PP\{\Ecal_1^c\}-\PP\{\Ecal_2^c\}-\PP\{\Ecal_4^c\}\ge 1-4d^{-\beta},
$$
and appropriate re-scaling of the constants.

To show (\ref{eq:bdlse2}), we first derive, via an identical argument to Step 2, that
\begin{eqnarray}
\nonumber
\|\Delta\|_{L_2(\Pi_n)}^2&\le& \|\Delta\|_{L_2(\Pi)}^2+C_{11}\zeta^{-{4\alpha\over 2\alpha-1}}\eta_qu^2+C_{11}\zeta^{{4\alpha\over 2\alpha+1}}n^{-(1-\max\{{q\over 2},{1\over 2\alpha+1}\})}\\
\nonumber
&&\qquad +C_{11}\sqrt{(\beta+1)\log d\over n}\left(\left({\log d\over n}\right)^{-q/4}u+2\left({\log d\over n}\right)^{1-q\over 2}\right)\\
\nonumber
&&\qquad +C_{11}n^{-1/2}e^{-d}\\
&&\qquad +C_{11}\left(u\sqrt{(\beta+1)\log d\over n}+{(\beta+1)\log d\over n}\right),
\label{eq:bdleftfinal2}
\end{eqnarray}
for some constant $C_{11}>0$. Together with (\ref{eq:bdlse}), this implies (\ref{eq:bdlse2}).


\section*{Appendix A -- Proof of Lemma \ref{le:emp1}}

An application of Talagrand's concentration inequality yields, with probability at least $1-e^{-t}$
$$
R_{jn}(u)\le 2\left(\EE R_{jn}(u)+u\sqrt{t\over n}+{t\over n}\right).
$$
It is well known that there exists a numerical constant $C_1>0$
$$
\EE R_{jn}(u)\le \left\{\EE \left[R_{jn}(u)\right]^{2}\right\}^{1/2}\le C_1n^{-1/2}u^{1-{1\over 2\alpha}}.
$$
See, e.g., Mendelson (2002) or Koltchinskii (2011). In other words, with probability at least $1-e^{-t}$,
$$
R_{jn}(u)\le C_2\left(n^{-1/2}u^{1-{1\over 2\alpha}}+u\sqrt{t\over n}+{t\over n}\right)
$$
for some numerical constant $C_2>0$. We now make this inequality uniform over $u\in [0,1]$ via a peeling argument.

In particular, with probability at least $1-\exp(-\beta\log d-2\log j)$ for some constant $\beta>0$,
\begin{eqnarray*}
\sup_{\substack{\|h\|_{\Hcal_1}\le 1\\ 2^{-j}\le \|h\|_{L_2(\Pi_j)}\le 2^{-j+1}}} \left|{1\over n}\sum_{i=1}^n \sigma_i h(x_{ij})\right|&\le& R_{jn}(2^{-j+1})\\
&\le&C_2n^{-1/2}\biggl[(2^{-j+1})^{1-{1\over 2\alpha}}+2^{-j+1}(\beta\log d+2\log j)^{1/2}\\
&&\qquad +n^{-1/2}(\beta\log d+2\log j)\biggr].
\end{eqnarray*}
By union bound, there exists a constant $C_3>0$ such that
$$
R_{jn}(u)\le C_3n^{-1/2}\left(u^{1-{1\over 2\alpha}}+u\sqrt{\beta\log d}+{\beta\log d\over \sqrt{n}}\right),
$$
holds for any $u\in (e^{-d(2\alpha/(2\alpha-1))},1]$, with probability at least
$$
1-\sum_{j=1}^{\lceil 2\alpha d\log_2e/(2\alpha-1)\rceil} \exp(-\beta\log d-2\log j)\ge 1-2d^{-\beta}.
$$
On the other hand, when $u\le e^{-d(2\alpha/(2\alpha-1))}$,
\begin{eqnarray*}
R_{jn}(u)&\le& R_{jn}(e^{-d(2\alpha/(2\alpha-1))})\\
&\le& C_2n^{-1/2}\left(e^{-d}+e^{-d(2\alpha/(2\alpha-1))}\sqrt{\beta\log d}+{\beta\log d\over \sqrt{n}}\right)\\
&\le& 2C_2n^{-1/2}\left(e^{-d}+{\beta\log d\over \sqrt{n}}\right),
\end{eqnarray*}
with probability at least $1-d^{-\beta}$, for sufficiently large $d$. In summary, there exists a constant $C_4>0$ such that
$$
R_{jn}(u)\le C_4n^{-1/2}\left(u^{1-{1\over 2\alpha}}+u\sqrt{\beta\log d}+{\beta\log d\over \sqrt{n}}+e^{-d}\right),
$$
uniformly over all $u\in [0,1]$ with probability at least $1-3d^{-\beta}$.

\section*{Appendix B -- Proof of Lemma \ref{le:emp}}

Note that
$$
\int_0^u \left[\log \Ncal(\Bcal_1(\Hcal_1), \delta, \|\cdot\|_{L_\infty})\right]^{1/2}du \le c_\alpha\delta^{1-{1\over 2\alpha}}.
$$
Therefore, there exist constants $C_1, C_2>0$ such that for any fixed $u\in [0,1]$
$$
\PP\left\{\hat{Z}_{jn}(u)\le C_1n^{-1/2}\left(u^{1-{1\over 2\alpha}}+ut^{1/2}\right)\right\}\le C_2\exp\left[-(u^{-1/\alpha}+t)\right].
$$
See, e.g., van de Geer (2000; Corollary 8.3). The rest of the proof follows a similar peeling argument as that for Lemma \ref{le:emp1} and is omitted for brevity. 

\end{document}